\documentclass[11pt,reqno]{amsart}

\usepackage[T1]{fontenc}
\usepackage{mathptmx,microtype,hyperref,graphicx,dsfont}
\usepackage{amsthm,amsmath,amssymb}
\usepackage[foot]{amsaddr}
\usepackage{enumitem}
\usepackage[nameinlink,capitalise]{cleveref}
\usepackage[font=footnotesize,margin=2cm]{caption}

\crefname{theorem}{Theorem}{Theorems}
\crefname{thm}{Theorem}{Theorems}
\crefname{mainthm}{Theorem}{Theorems}
\crefname{lemma}{Lemma}{Lemmas}
\crefname{lem}{Lemma}{Lemmas}
\crefname{remark}{Remark}{Remarks}
\crefname{claim}{Claim}{Claims}
\crefname{subclaim}{Sub-claim}{Sub-claims}
\crefname{prop}{Proposition}{Propositions}
\crefname{proposition}{Proposition}{Propositions}
\crefname{defn}{Definition}{Definitions}
\crefname{corollary}{Corollary}{Corollaries}
\crefname{conjecture}{Conjecture}{Conjectures}
\crefname{question}{Question}{Questions}
\crefname{chapter}{Chapter}{Chapters}
\crefname{section}{Section}{Sections}
\crefname{figure}{Figure}{Figures}
\newcommand{\crefrangeconjunction}{--}

\theoremstyle{plain}

\newtheorem{thm}{Theorem}
\newtheorem*{thm*}{Theorem}
\newtheorem{lemma}[thm]{Lemma}

\newtheorem{prop}[thm]{Proposition}

\theoremstyle{definition}
\newtheorem{defn}[thm]{Definition}

\theoremstyle{remark}

\newcommand{\eps}{\varepsilon}
\renewcommand{\P}{{\bf P}}

\newcommand{\p}{\sqrt{\alpha/(n\log{n})}}

\renewcommand{\l}{\langle}
\renewcommand{\r}{\rangle}

\newcommand{\G}{{\mathcal G}}
\newcommand{\Gnp}{\G_{n,p}}

\newcommand{\CP}{{\mathcal C}}

\author[B. Kolesnik]{Brett Kolesnik}

\address{Department of Statistics, University of California, Berkeley}

\email{bkolesnik@berkeley.edu}

\begin{document}

\title{Sharp threshold for $K_4$-percolation}

\begin{abstract}
We locate the critical threshold
$p_c$
at which  
it becomes likely that the complete graph $K_n$ 
can be obtained from the 
Erd{\H{o}}s--R{\'e}nyi graph $\G_{n,p}$
by iteratively completing copies of $K_4$ minus an 
edge. 
This refines work of Balogh, Bollob{\'a}s and Morris  
that bounds the threshold up to multiplicative
constants. 
\end{abstract}

\maketitle


\section{Introduction}\label{S_intro}

Triangles play an important role in networks.
For instance, the concept of triadic closure 
\cite{S08,G73} from social network theory 
is the observation that if there are edges (e.g., representing friendship)
between vertices $x,y$ and $x,z$, 
then the edge $y,z$ 
(if not already present) 
 is likely to be added eventually. 
This gives rise
to the special case $H=K_3$ in the process called 
{\it $H$-graph bootstrap percolation}
introduced by 
Bollob\'as~\cite{B67} (under the name of {\it weak saturation}). 
Let $\l G\r_H$ denote the graph obtained from $G$ by 
iteratively
completing 
copies of $H$ minus an edge. 
A graph $G$ is said to {\it $H$-percolate}
if all missing edges are eventually added, 
that is, if $\l G\r_H$ is the complete graph on the vertices
of $G$. 

Following  
Balogh, Bollob{\'a}s and Morris~\cite{BBM12}, 
we suppose that the underlying network is 
the 
Erd{\H{o}}s--R{\'e}nyi~\cite{ER59} graph,
that is,  the random subgraph $\G_{n,p}$ of the complete
graph $K_n$ where edges are included independently with probability $p$.
The  critical threshold,
at which  
$\Gnp$ is likely to $H$-percolate,
is defined formally as 
\[
  p_c(n,H) 
  = \inf\left\{p>0 : \P(\l\Gnp\r_H = K_n) \ge 1/2 \right\}.
\]
A graph $K_3$-percolates if and only if it is connected,
so this case follows
by standard results \cite{ER59}. 
In this work, we focus on the next case, $H=K_4$. 
This is a more stringent version of 
triadic closure, where edges $u,v$ are added only if $u$
and $v$ are incident to triangles that share an edge
(e.g., people become friends if they have mutual friends
who are friends). 
In \cite{BBM12},  
$p_c(n,K_4)$ is estimated up to multiplicative constants. 
Our main result locates the sharp threshold. 

\begin{thm}\label{T_K4}
$p_c(n,K_4)\sim1/\sqrt{3n\log n}$. 
\end{thm}

\subsection{Outline}\label{S_out}
The upper bound is proved in \cite{AK18}, via  
a connection with classical {\it $2$-neighbor bootstrap percolation}
\cite{PRK75,CLR79,JLTV12}, which we now explain.  
In this model, vertices 
are infected if they have at least 2 infected neighbors. 
Suppose that some set $I$ of vertices in a graph $G=(V,E)$
are initially infected. Let $\l I,G\r_2$ denote the 
set of eventually infected vertices. 
If all vertices are eventually infected, $\l I,G\r_2=V$, 
we say that $I$ is {\it contagious} for $G$. It is easy to see 
(by induction) 
that 
if some edge in $G$ is contagious, then $G$ 
will $K_4$-percolate. 
Therefore, the upper bound in \cref{T_K4}
follows since, as shown in  \cite{AK18}, $1/\sqrt{3n\log n}$
is the sharp threshold for the existence of such a 
{\it seed edge} in $\Gnp$.

To prove the lower bound in \cref{T_K4}, we essentially 
show that none of the other ways in which $\Gnp$
can percolate are more likely. The analysis is 
somewhat involved, as there are many ways in which 
percolating subgraphs of $\Gnp$ can ``merge'' to form
larger percolating subgraphs. 
Similar issues are involved, for instance, 
 with the pioneering work of Holroyd \cite{H03}.

The key to overcoming this difficulty, in the current work, 
is the observation that if a graph $G$ percolates, 
then the subgraph $C$ obtained 
by successively deleting vertices of degree 
2 also percolates. 
We call $C$ the {\it core} of $G$. The case that $C$ is 
a  seed edge
is described above. Otherwise, $C$ has minimum degree
at least 3, in which case we call $C$ a {\it $3$-core}. Hence, 
a percolating graph $G$ is either a {\it seed graph,}
or else it has a $3$-core. In either case, the vertex set 
$V(C)$ is contagious for $G$. 

There are two other main ingredients in the 
proof of the lower bound. 
First, 
by a detailed combinatorial analysis, 
based on the clique process 
(see \cref{S_CP} below)
defined in \cite{BBM12}, we show 
that there are at most $(2/e)^q q!q^q$ 
percolating 
3-cores of size $q$. 
Then, with this at hand, we utilize the following tail estimates
\cite{AK17a} (which complement the central limit theorems in 
\cite{JLTV12}).
 
Let $P(q,k)$ denote the probability 
that for a given set $I\subset[n]$ (independent of $\Gnp$), with $|I|=q$, 
we have that 
$|\l I,\Gnp\r_2|\ge k$.

\begin{lemma}[{\cite{AK17a}}]
\label{L_P}
Fix $\alpha>0$ and put  
$p=\sqrt{\alpha/(n\log{n})}$. Let $\eps\in[0,1)$ and 
$\beta\in[\beta_\eps,1/\alpha]$, where  
$\beta_\eps=(1-\sqrt{1-\eps})/\alpha$. 
Put $k_\alpha=\alpha^{-1}\log{n}$
and $q_\alpha=(2 \alpha)^{-1}\log{n}$. 
Suppose that $q/q_\alpha\to\eps$
and $k/k_\alpha\to\alpha\beta$ as $n\to\infty$. 
Then
$
P(q,k)
= 
n^{\xi+o(1)}
$, 
where 
\[
\xi=-\frac{\alpha\beta^2}{2}+
\begin{cases}
(2\alpha\beta-\eps)(2\alpha)^{-1}
\log(
e(\alpha\beta)^2/(2\alpha\beta-\eps)
)
& \beta\in[\beta_\eps,\eps/\alpha) \\
\beta\log(\alpha\beta)
-\eps(2\alpha)^{-1}\log(\eps/e)
& \beta\in[\eps/\alpha,1/\alpha].
\end{cases}
\]
\end{lemma}

(This follows by the main result in \cite{AK17a}, 
setting $r=2$ and 
replacing the parameters $\vartheta,\alpha,\beta$
therein 
with $k_\alpha,\eps,\alpha\beta$, respectively.) 

Using this, together with the upper bound 
$(2/e)^q q!q^q$ 
for percolating 3-cores
of size $q$, we  
argue (see \cref{S_smallC}) that, when $p$ is sub-critical,  
the expected number of percolating subgraphs 
of $\Gnp$ of size $k=\beta\log{n}$,  
for $\beta\in[\beta_\eps,1/\alpha]$, 
with a core of size $q\le (3/2)\log n$ 
is bounded by $n^{\mu+o(1)}$, where
\[
\mu(\alpha,\beta)=3/2
+\beta\log(\alpha\beta)-\alpha\beta^2/2.  
\] 
The almost sure non-existence of percolating 3-cores
of size $q\ge (3/2)\log n$ in $\Gnp$ is handled separately (see \cref{S_largeC}), 
by showing that such a graph would have to be created through a 
highly unlikely merging of other graphs of ``macroscopic'' size. 
This leads to the following result, yielding the lower
bound in \cref{T_K4}. 

\begin{thm}\label{T_beta*}
Fix $\alpha\in(0,1/3)$ and put 
$p=\sqrt{\alpha/(n\log{n})}$.
With high probability the largest cliques in $\l \Gnp \r_{K_4}$
are of size $(\beta_*+o(1))\log{n}$, 
where    
$\mu(\alpha,\beta_*)=0$.
\end{thm}

\subsection{Discussion}

The critical window for the
connectivity of $\Gnp$ is well-understood. 
With high probability 
$\Gnp$ is connected (hence $K_3$-percolating) 
if and only if it has no isolated vertices. 
If 
$p=(\log{n}+\eps)/n$,  
$\Gnp$ will $K_3$-percolate
with probability $\exp(-e^{-\eps})(1+o(1))$,
as $n\to\infty$. It would be interesting to obtain
similarly detailed information for $K_4$-percolation.

Estimates for $p_c(n,K_r)$ up to poly-logarithmic factors
are obtained in \cite{BBM12}. 
More recently, the threshold $p_c$   
has been located up to constant
factors \cite{BK20}. Interestingly, 
the connection with classical bootstrap percolation 
described above does not lead to the critical threshold
when $r\ge5$. Instead, near $p_c$, $\Gnp$ percolates
in some other way, that is still not fully understood. 

Although $H$-percolation
can in general   behave quite differently than the present case $H=K_4$, 
we think the ideas in this work will be useful in 
improving the bounds for $p_c$ in other cases
of interest.

\section{The clique process}\label{S_CP}

The \emph{clique process} \cite{BBM12} 
describes the $K_4$-percolation dynamics
in a way that is amenable to analysis. 

\begin{defn}
Three graphs $G_i=(V_i,E_i)$ 
\emph{form a triangle}
if there are distinct vertices $x,y,z$ such that 
$x\in V_1\cap V_2$, $y\in V_1\cap V_3$ and 
$z\in V_2\cap V_3$. If $|V_i\cap V_j|=1$ for all $i\neq j$, 
we say that they \emph{form exactly one triangle}. 
\end{defn}

In \cite{BBM12} the following observation is made.

\begin{lemma}\label{L_merge}
Suppose that $G_i=(V_i,E_i)$
percolate. 
\begin{enumerate}[nolistsep,label=(\roman*)]
\item  If the $G_i$ form a triangle then 
$G_1\cup G_2\cup G_3$ percolates.
\item If $|V_1\cap V_2|\ge2$ then $G_1\cup G_2$
percolates. 
\end{enumerate}
\end{lemma}

This leads to the following process. 

\begin{defn}
A \emph{clique process} for a graph $G$
is a sequence $(\CP_t)_{t=0}^\tau$ of sets
of subgraphs of $G$ such that: 
\begin{enumerate}[nolistsep,label={(\roman*)}]
\item $\CP_0=E(G)$ is the edge set of $G$.  
\item For each $t<\tau$, $\CP_{t+1}$ is obtained 
from $\CP_t$ by either 
(a) merging two subgraphs $G_1,G_2\in\CP_t$
with at least two common vertices, or 
(b) merging three
subgraphs 
$G_1,G_2,G_3\in\CP_t$
that form exactly one triangle. 
\item $\CP_\tau$ is such that no further operations 
as in (ii) are possible. 
\end{enumerate}
\end{defn}

The reason for the name is that (by induction), for any $t\le \tau$ and $H\in\CP_t$, 
$H$ percolates, and 
hence 
$\l H\r_{K_4}$ is 
a clique in $\l G\r_{K_4}$. 

The description above is slightly modified 
from that presented in \cite{BBM12}, as we note that 
if three percolating graphs form more than one triangle, then they can be 
merged by applying  \cref{L_merge}(ii) twice. 
Therefore, for convenience, we reserve the use of 
\cref{L_merge}(i) in a clique process for the case that 
{\it exactly} one triangle is formed.
This simplifies the combinatorial analysis 
in \cref{S_core+deg2s} below.

Finally, let us record the following observation, 
see \cite{BBM12}.

\begin{lemma}\label{L_CP}
Let $G$ be a finite graph and 
$(\CP_t)_{t=0}^\tau$ a clique process for $G$. 
For each $t\le \tau$, $\CP_t$ is a set of edge-disjoint, 
percolating subgraphs of $G$. Furthermore, 
$\l G\r_{K_4}$ is the edge-disjoint, 
triangle-free union of cliques 
$\bigcup_{H\in \CP_\tau}\l H\r_{K_4}$. 
Hence $G$ percolates if and only if $\CP_\tau=\{G\}$. 
In particular, $\CP_\tau=\CP'_{\tau'}$
for any two clique processes 
$(\CP_t)_{t=0}^\tau$ and $(\CP'_t)_{t=0}^{\tau'}$
for $G$. 
\end{lemma}

\subsection{Consequences}\label{S_CPcons}

The following    
consequences of \cref{L_CP}, derived in \cite{BBM12}
using the clique process,
play a crucial role in the current work.  

\begin{lemma}\label{L_Emin}
If $G=(V,E)$ percolates then 
$|E|\ge2|V|-3$. 
\end{lemma}

\begin{defn}
We call $|E|-(2|V|-3)$ the \emph{excess} 
of a graph 
$G=(V,E)$. 
A  graph is 
{\it edge-minimal} if its excess is 0. 
\end{defn}

To prove \cref{L_Emin}, the following observations
are made in \cite{BBM12}. 

\begin{lemma}\label{L_Emerge}
Suppose that  $G_i=(V_i,E_i)$
percolate. 
\begin{enumerate}[nolistsep,label={(\roman*)}]
\item If the $G_i$ form exactly one triangle, 
then the excess of $G_1\cup G_2\cup G_3$ 
is the sum of the excesses of the $G_i$. 
\item If $|V_1\cap V_2|=m\ge2$, then the excess of 
$G_1\cup G_2$ is the sum of the excesses of the $G_i$, 
plus $2m-3>0$.  
\end{enumerate}
\end{lemma}

Hence, if $G$ is an edge-minimal percolating graph,
then every step of any clique process for $G$ involves
merging three subgraphs that form exactly one triangle. 
The simplest example of this is when two of the $G_i$
are a single edge sharing a common vertex. 
If all steps of a clique process for $G$ are of this form, then
$G$ is a seed graph, as defined in \cref{S_out} above.

Finally, since at most three subgraphs are merged
in any step of a clique process,  
we have the following 
Aizenman--Lebowitz \cite{AL88}
type condition.

\begin{lemma}\label{L_k3k}
Let $G$ be a graph  
and $k\ge1$. If $G$ has  
no percolating subgraphs 
of size $k'\in[k,3k]$ then $G$ has no
percolating subgraphs larger than $k$. 
\end{lemma}

\section{Combinatorial bounds}\label{S_core+deg2s}

We first address the issue of estimating the 
number of percolating graphs with various structural properties. 
Most crucially, we require reasonably  sharp estimates
for the number of percolating graphs with few vertices of 
degree 2. 
The proofs of the main results in this section 
\cref{L_Iell,L_Iellq} are fairly straightforward, but rather
laborious. As such, we only sketch
the main ideas in the proofs in this section. 
The proofs appear in \cref{A_Iell,A_Iellq}
below. 

\begin{defn}
A percolating graph $G$ is \emph{irreducible}
if removing any edge from $G$ results in a 
non-percolating graph. 
\end{defn}

Note that a graph can be irreducible, but not edge-minimal. 

Clearly, a graph $G$ percolates if and only if
it has an irreducible percolating subgraph 
$G'\subset G$ such that $V(G)=V(G')$.

Next, we observe that if a vertex of degree 2 is 
removed from a percolating graph, the resulting subgraph 
still percolates. 
This follows by arguments in \cite{B67}, however, 
this article is not widely accessible. 
For completeness, a proof using the clique process
is given in \cref{A_deg2s}
below. 

\begin{lemma}\label{L_deg2s}
Suppose that $G$ percolates
and $v\in V(G)$ is of degree $2$. Then the subgraph $G_v\subset G$
induced by $V-\{v\}$
percolates.
\end{lemma}

As discussed in \cref{S_out},  
$\l I,G\r_2$ denotes the set of vertices eventually infected 
by the 2-neighbor dynamics on $G$, when 
$I$ is initially infected. 

\begin{defn}
Similarly, 
for a subgraph $H\subset G$, 
we write $\l H,G\r_2$ to denote
the subgraph of $G$ induced by $\l V(H),G\r_2$.
\end{defn}

By \cref{L_merge}(i) and induction, if $H\subset G$ percolates, 
then so does $\l H,G\r_2$.

The following is an immediate consequence
of \cref{L_deg2s,L_Emerge}. 

\begin{lemma}\label{L_core+deg2s}
Let $G$ be an irreducible percolating graph. 
Then either: 
\begin{enumerate}[nolistsep,label={(\roman*)}]
\item $G=\l e,G\r_2$ for some edge $e\in E(G)$, 
or else, 
\item $G= \l C,G\r_2$ for some percolating  $C\subset G$
of minimum degree at least $3$.
\end{enumerate}
Furthermore: 
\begin{enumerate}[nolistsep,label={(\roman*)}]
\item[(iii)] the excess of $G$ is equal to the excess of $C$. 
\end{enumerate}
\end{lemma}

In the first case, 
$G$ is a seed graph and $e$ is a seed edge.
Such a graph (if irreducible) is 
edge-minimal. 
In the latter case, $C$ is the $3$-core of $G$. 
If $G=C$, we say that 
$G$ is a $3$-core.

\begin{figure}[h]
\includegraphics[scale=0.8]{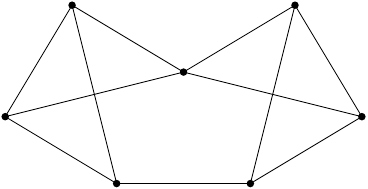}
\caption{
The smallest irreducible percolating $3$-core.
}
\label{F_3core}
\end{figure}

It is easy to see that all irreducible 
percolating graphs 
on $2<k\le 6$
vertices have a vertex of degree $2$.
There are, however, edge-minimal percolating graphs 
of size $k=7$ (and larger) with no vertices of degree $2$, 
e.g., see \cref{F_3core}.

\subsection{Basic estimates}\label{S_est}

In this section, we 
use \cref{L_core+deg2s} to
obtain upper bounds for the number of 
irreducible percolating graphs.
For such a graph $G$, 
the relevant quantities are  its size, 
the number of vertices in $G$ of degree 2,  
the size of its core $C\subset G$, 
and its number of excess edges.

\begin{defn}
Let $I^\ell_q(k,i)$ be the number of  
labelled, irreducible percolating graphs $G$
of size $k$ with an excess of $\ell$ edges, $i$ vertices of 
degree $2$, and a core $C\subset G$ of 
size $q$. 
If $i=0$, and hence $q=k$, 
we simply write $C^\ell (k) = I^\ell_k(k,0)$. 
If $\ell=0$, we write
$I_q(k,i)$ and $C(k)$.  
\end{defn}

Note that $I_2(k,i)$ is the number of 
labelled, irreducible (and edge-minimal) seed graphs 
of size $k$ with $i$ vertices of degree $2$.

By \cref{L_core+deg2s}(iii), if a graph $G$ contributes
to $I^\ell_q(k,i)$ then its core  
has an excess of $\ell$ edges. 
As noted above, 
there are no irreducible $3$-cores
on $q\le 6$ vertices. Hence   
$I_q^\ell(k,i)=0$ if $2<q\le 6$.

\begin{defn}
We let $I^\ell(k,i)=\sum_{q} I_q^\ell(k,i)$
denote the number of  
labelled, irreducible graphs $G$
of size $k$, with an excess of $\ell$ edges
and  $i$ vertices of 
degree $2$.
If $\ell=0$, we  write $I(k,i)$.
\end{defn}

We obtain the following 
estimates for $I^\ell(k,i)$, 
assuming the excess is 
$\ell\le3$.  The method of proof could presumably 
(with additional work) provide 
bounds for larger $\ell$, however, fortunately, 
percolating graphs with 
a larger excess can be 
dealt with using less accurate estimates
(see \cref{L_extraE} below).    

\begin{lemma}\label{L_Iell}
For all $k\ge2$, $\ell\le3$ and relevant $i$,
we have that 
\[
I^\ell(k,i)\le (2/e)^k k!k^{k+2\ell+i}.
\]
In particular, 
$C^\ell(k)\le(2/e)^k k!k^{k+2\ell}$.
\end{lemma}

Note that, for small values of $i$, 
this is much smaller than the total number of 
seed graphs of size $k$, which
in \cite{AK18} is shown to be  
roughly equal to $k!k^k$.

See \cref{A_Iell} below for the proof.  
The argument is 
quite 
lengthy, as 
there are several cases 
(increasing in $\ell$) 
to consider, 
depending on the nature
of the last step in the clique process. 
Before moving on, we sketch the main ideas. 

First, we note that the cases $i>0$ follow by a simple induction, 
since if $G$ has $i$ vertices of degree 2, 
then removing such a vertex from $G$
results in a graph with $j\in\{i,i\pm1\}$ 
vertices of degree 2. 
Analyzing this case leads to the optimal constant $2/e$. 
The case of 3-cores $i=0$
is the heart
of the proof.  
The following observation is the key: 
If $G$ is a percolating $3$-core, 
then in the last step of a clique process, 
either (i) three graphs $G_1,G_2,G_3$ are merged that
form exactly one triangle on $T=\{v_1,v_2,v_3\}$, 
or else (ii)
two graphs $G_1,G_2$ are merged 
that share $m\ge2$
vertices $S=\{v_1,v_2,\ldots,v_m\}$.  
If 
some $G_j$ has a vertex $v$ of degree 2, 
then necessarily $v\in T$ in case (i), 
or $v\in S$ in case (ii) 
(as else  $G$ would have a vertex of degree 2). 
In other words, if a percolating $3$-core is formed
by merging graphs with vertices
of degree 2, then all such vertices belong to the 
triangle that they form or the set of vertices that they share. 
These observations provide enough control on the combinatorics
to allow for an inductive proof of the bounds in 
\cref{L_Iell}.

\subsection{Sharper estimates}\label{S_est-q}

Next, using \cref{L_Iell} as a starting point, 
we obtain the following upper bounds for $I_q^\ell(k,i)$.  

\begin{lemma}\label{L_Iellq}
Fix $\eps>0$. 
For some constant $\vartheta(\eps)\ge1$, the following holds. 
For all $k\ge2$, $\ell\le3$,  
and relevant $q,i$,
we have that  
\[
I^\ell_q(k,i)\le \vartheta \psi_\eps(q/k)^k k!k^{k+2\ell+i}
\]
where
\[
\psi_\eps(y)=
\max
\{
3/(2e)+\eps,(e/2)^{1-2y}y^2
\}.
\] 
\end{lemma}

This lemma 
improves upon \cref{L_Iell} only 
when $\eps<1/(2e)$, as otherwise 
$\psi_\eps(y)\ge2/e$ for all $y$. 
On the other hand, for any given $\eps<1/(2e)$, we have that $\psi_\eps(y)$ is 
non-decreasing and  $\psi_\eps(y)\to 2/e$ as $y\uparrow 1$.
Note that 
$\psi_\eps(y)=3/(2e)+\eps$
for $y\le y_*$ and $\psi_\eps(y)=(e/2)^{1-2y}y^2$
for $y>y_*$, where 
\begin{equation}\label{E_y*}
3/(2e)+\eps=(e/2)^{1-2y_*}y_*^2.
\end{equation}
We define $ y_0=y_*(0)\approx0.819$, and note that 
$y_*(\eps)\downarrow y_0$, as $\eps\downarrow0$.

The main ideas in the proof are as follows: 
First, we note that the case $i=k-q$ 
follows essentially directly by \cref{L_Iell}. 
We establish the remaining cases by induction, 
noting that if a graph $G$ contributes to 
$I^\ell_q(k,i)$ and $i<k-q$, then 
there is a vertex $v$ in $G$ of degree 2 
with a neighbor that is not in the core $C\subset G$.
Therefore, either 
(i) some neighbor of $v$ is of degree 2
in $G_v$, or else 
(ii) there are vertices $u\neq w$ of degree 2 in $G$ with a common 
neighbor that is not in $C$. 
Beyond these observations, the proof 
is mostly calculus, see \cref{A_Iellq}
below.

\section{Proof of \cref{T_beta*}}\label{S_proof}

With our key 
\cref{L_P,L_Iell,L_Iellq} at hand, 
we turn to the proof of
\cref{T_beta*}.
The argument is divided into two parts
\cref{S_smallC,S_largeC} where, respectively, 
percolating subgraphs of $\Gnp$ with small 
and large cores are considered.

\subsection{Percolating subgraphs 
with small cores}
\label{S_smallC}

First, we show that for sub-critical $p$,
with high probability  
$\Gnp$ has no percolating subgraphs 
significantly larger that $\beta_*\log{n}$ 
with a small core.

\begin{prop}\label{P_smallC}
Fix $\alpha\in(0,1/3)$ and put $p=\p$. 
Then, for any $\delta>0$, 
with high probability $\Gnp$ has no
irreducible percolating subgraphs $G$ of size 
$k=\beta\log{n}$,  
for $\beta\ge\beta_*+\delta$, 
 with a core $C\subset G$
of size $q\le (3/2)\log{n}$. 
\end{prop}

First, we note  
that $\beta_*$
in \cref{T_beta*} is well-defined. 

\begin{lemma}\label{L_beta*}
Fix $\alpha\in(0,1/3)$. For $\beta>0$, let 
  \[
\mu(\alpha,\beta)
=3/2
+\beta\log(\alpha\beta)-\alpha\beta^2/2.
\] 
 The function $\mu(\alpha,\beta)$ is decreasing 
 in $\beta$, with a unique zero 
 $\beta_*\in(0,3)$. 
\end{lemma}

\begin{proof}
Differentiating $\mu(\alpha,\beta)$ 
with respect to $\beta$, we obtain  
$1+\log(\alpha\beta)-\alpha\beta$.  
Since $\log{x}<x-1$ for all positive $x\neq1$,
we find that $\mu(\alpha,\beta)$ is decreasing in $\beta$. 
Moreover, since $\alpha<1/3$, we have that  
$\mu(\alpha,3)<(3/2)(3\alpha-1)<0$.
The result follows, noting that 
$\mu(\alpha,\beta)\to3/2>0$ as $\beta\downarrow0$. 
 \end{proof}

Recall that the bounds in 
\cref{L_Iell,L_Iellq}  apply only to graphs with 
an excess of $\ell\le3$ edges. 
For graphs with larger excess, we will 
apply the following result. 

\begin{lemma}\label{L_extraE}
Fix $\alpha\in(0,1/3)$ and put $p=\p$. 
Then
with high probability $\Gnp$ contains no subgraph 
of size $k=\beta\log{n}$ with an excess of $\ell$ edges, 
for any $\beta\in(0,2]$ and $\ell>3$, 
or any $\beta\in(0,9]$ and $\ell>27$.  
\end{lemma}

\begin{proof}
The expected number of subgraphs of size 
$k=\beta\log{n}$ in $\Gnp$ 
with an excess of $\ell$ edges is bounded by
\[
{n\choose k}{{k\choose2}\choose 2k-3+\ell}p^{2k-3+\ell}
\le 
\left(\frac{e^3}{16}knp^2\right)^k
\left(\frac{e}{4}kp\right)^{\ell-3}\le n^\nu\log^{\ell}n
\]
where 
\[
\nu(\beta,\ell)=-(\ell-3)/2
+\beta\log(\alpha\beta e^3/16).
\]
Note that $\nu$ is convex in $\beta$ and 
$\nu(\beta,\ell)\to-(\ell-3)/2$ as $\beta\downarrow0$. 
Note also that
\[
2\log(2/3\cdot e^3/16)\approx-0.356<0
\]
and 
\[
9\log(9/3\cdot e^3/16)\approx 11.934<12.
\]
Therefore, since $\alpha<1/3$, 
$\nu(2,\ell)<-(\ell-3)/2$ and 
$\nu(9,\ell)<-(\ell-27)/2$. 
The result follows. 
 \end{proof}

\begin{defn}
Let 
$E(q,k)$ denote the 
expected number of irreducible 
percolating cores  
$C\subset \Gnp$ of size $q$  
such that $|\l C,\Gnp\r_2|\ge k$.
\end{defn}

Combining \cref{L_P,L_Iell,L_extraE}, 
we obtain the following 
estimate. Recall
$\beta_\eps,k_\alpha,q_\alpha$ as in \cref{L_P}, 
and $\mu$ in \cref{L_beta*}. 

\begin{lemma}\label{L_mueps}
Fix $\alpha\in(0,1/3)$ and put $p=\p$. 
Let $\eps\in[0,3\alpha]$ and
 $\beta\in[\beta_\eps,1/\alpha]$. 
Suppose that $q/q_\alpha\to\eps$
and $k/k_\alpha\to\alpha\beta$ as $n\to\infty$. 
Then 
$
E(q,k)\le n^{\mu_\eps+o(1)}
$, 
where 
$\mu_\eps(\alpha,\beta)=\mu(\alpha,\beta)$ for 
$\beta\in[\eps/\alpha,1/\alpha]$, 
\[
\mu_\eps(\alpha,\beta)
= 
\mu(\alpha,\beta)
-\beta\log(\alpha\beta)
+
\frac{\eps}{2\alpha}\log(\eps/e)
+
\frac{2\alpha\beta-\eps}{2\alpha}\log\left(
\frac{e(\alpha\beta)^2}{2\alpha\beta-\eps}
\right)
\]
for $\beta\in[\beta_\eps,\eps/\alpha]$. 
\end{lemma}

\begin{proof}
By \cref{L_extraE}, it suffices to 
show that, for all $\ell\le3$, 
we have that $E^\ell(q,k)\le n^{\mu_\eps+o(1)}$, where
$E^\ell(q,k)$ is the expected number of irreducible percolating
cores $C\subset\Gnp$ of size $q=\eps(2\alpha)^{-1}\log{n}$
with an excess of $\ell$ edges, such that 
$|\l C,\Gnp\r_2|\ge k=\beta\log{n}$. 
For such $\ell$, by
\cref{L_P,L_Iell}, 
we find that 
\begin{align*}
E^\ell(q,k)
&\le
{n\choose q} C^\ell(q) p^{2q-3+\ell} P(q,k)\\
&\le 
q^{2\ell}p^{\ell-3}\left(\frac{2}{e}qnp^2\right)^{q}
P(q,k)\le n^{\nu+o(1)}
\end{align*}
where
\[
\nu = 
3/2
+\eps(2\alpha)^{-1}\log(\eps/e)
+\xi_\eps(\alpha,\beta)=\mu_\eps(\alpha,\beta).
\qedhere
\]
 \end{proof}

Having established \cref{L_mueps}, 
we aim to prove \cref{P_smallC}
by the first moment method. 
We first show that for some $\eps_*\in(0,3\alpha)$,
with high probability there are no irreducible percolating 
cores in $\Gnp$ of size $\eps(2\alpha)^{-1}\log{n}$, 
with $\eps\in(\eps_*,3\alpha]$. 
We record a slightly more general result, allowing 
for $i=O(1)$ vertices of degree 2, 
as this will be required in \cref{S_largeC} below. 

\begin{lemma}\label{L_eps*}
Fix $\alpha\in(0,1/3)$ and put $p=\p$. 
Fix some $i_*\ge0$. 
Define $\eps_*\in(0,3\alpha)$ implicitly by   
$3/2+\eps_*(2\alpha)^{-1}\log(\eps_*/e)=0$. Then, for any $\eta>0$, 
with high probability $\Gnp$ has no irreducible percolating 
subgraphs $G$ of size $k=\eps(2\alpha)^{-1}\log{n}$
with $i\le i_*$ vertices of degree $2$, for  
$\eps\in[\eps_*+\eta,3\alpha]$. 
\end{lemma}

\begin{proof}
By \cref{L_extraE}, it suffices to consider subgraphs $G$
with excess  $\ell\le 3$.
By \cref{L_Iell}, the expected number
of such subgraphs is bounded by 
\[
{n\choose k}p^{2k-3+\ell} I^\ell(k,i)
\le k^{2\ell+i}p^{\ell-3} \left(\frac{2}{e}knp^2\right)^k
\le n^{\nu+o(1)}
\]
where
$
\nu(\eps) 
= 
3/2+\eps(2\alpha)^{-1}\log(\eps/e)
$.
Since $\nu$ 
is decreasing in $\eps<1$, 
$\nu\to3/2>0$ as $\eps\downarrow0$, 
and 
$
\nu(3\alpha)=(3/2)\log(3\alpha)<0
$, 
the lemma follows. 
 \end{proof}

Next, we use \cref{L_mueps} to rule out the remaining cases
$\eps\le\eps_*+\eta$ (where $\eta>0$ is a small constant, 
to be determined
below). In order to apply \cref{L_mueps}, we first verify that,  
for such $\eps$, 
we have 
$\beta_*\ge\beta_\eps$. 

\begin{lemma}\label{L_betaeps}
Fix $\alpha\in(0,1/3)$. Let 
$\beta_\eps,\beta_*,\eps_*$
be as 
in \cref{L_P,L_beta*,L_eps*}. Then, 
for some sufficiently small $\eta(\alpha)>0$, 
we have that $\beta_*\ge\beta_\eps$
for all $\eps\in[0,\eps_*+\eta]$. 
\end{lemma}

\begin{proof}
By \cref{L_beta*} and the continuity of $\mu(\alpha,\beta_\eps)$
in $\eps$, it suffices to show that 
$\mu(\alpha,\beta_\eps)>0$, for all $\eps\in[0,\eps_*]$. 
Let $\delta_\eps=1-\sqrt{1-\eps}$, so that 
$\beta_\eps=\delta_\eps/\alpha$. 
Note that 
\[
\mu(\alpha,\beta_\eps) 
=
3/2
+(2\alpha)^{-1}(2\delta_\eps\log\delta_\eps-\delta_\eps^2). 
\] 
Therefore, by the bound $\log{x}\le x-1$, 
\[
\frac{\partial}{\partial \eps}\mu(\alpha,\beta_\eps) 
=(2\alpha)^{-1}
(1+\log(\delta_\eps)/(1-\delta_\eps))
\le 0. 
\]
It thus suffices to verify that $\mu(\alpha,\beta_{\eps_*})>0$. 
To this end note that, by the definition of $\eps_*$
(see \cref{L_eps*}),
\[
\mu(\alpha,\beta_{\eps_*})
=
(2\alpha)^{-1}(
2\delta_{\eps_*}\log\delta_{\eps_*}-\delta_{\eps_*}^2
-\eps_*\log(\eps_*/e)). 
\]
By \cref{L_eps*}, we have that 
$\eps_*=\delta_{\eps_*}(2-\delta_{\eps_*})\in(0,1)$, and so
$\delta_{\eps_*}\in(0,1)$. Hence the lemma follows 
if we show that $\nu(\delta)>0$ for all $\delta\in(0,1)$, where
\[
\nu(\delta) 
= 
2\delta\log\delta-\delta^2
-\delta(2-\delta)\log(\delta(2-\delta)/e). 
\]
Note that
\[
\nu(\delta) /\delta
= 
\delta\log{\delta} - (2-\delta)\log(2-\delta) + 2(1-\delta). 
\]
Differentiating this expression with respect to $\delta$,
we obtain $\log(\delta(2-\delta))<0$, for all $\delta<1$. 
Noting that  $\nu(1)=0$,  
the lemma follows. 
 \end{proof}

It can be seen that, for all sufficiently large $\eps<\eps_*$, 
we have that $\beta_*<\eps/\alpha$, where $\mu_\eps\neq \mu$. 
Therefore, 
we require the following bound. 

\begin{lemma}\label{L_mueps<mu}
Fix $\alpha\in(0,1/3)$. Let $\eps\in[0,1)$ and 
$\beta_\eps,\mu_\eps$
be as 
in \cref{L_P,L_mueps}. 
Then $\mu_\eps(\alpha,\beta)\le\mu(\alpha,\beta)$, 
for all $\beta\in[\beta_\eps,1/\alpha]$.
\end{lemma}

\begin{proof}
Since $\mu(\alpha,\beta)=\mu_\eps(\alpha,\beta)$ 
for $\beta\in[\eps/\alpha,1/\alpha]$, we may assume
that $\beta<\eps/\alpha$. 
Let $\delta=\alpha\beta$. Then 
\[
\alpha(
\mu(\alpha,\beta) - 
\mu_\eps(\alpha,\beta))
=
\delta\log\delta
-
\frac{\eps}{2}\log(\eps/e)
-
\frac{2\delta-\eps}{2}\log\left(
\frac{e\delta^2}{2\delta-\eps}
\right). 
\]
Differentiating this expression with respect to $\delta$,
we obtain 
\[
\eps/\delta-1-\log(\delta/(2\delta-\eps))
\le0, 
\]
by the inequality $\log{x}\ge (x-1)/x$. 
Since $\mu(\alpha,\eps/\alpha)=\mu_\eps(\alpha,\eps/\alpha)$, 
the lemma follows. 
 \end{proof}

Finally, we prove the main result of this section. 

\begin{proof}[Proof of \cref{P_smallC}]
Let $\delta>0$ be given. By \cref{L_beta*}, we may assume 
that $\beta_*+\delta<1/\alpha$. 
If $\Gnp$ has an irreducible percolating  
subgraph $G$ of size $k\ge(\beta_*+\delta)\log{n}$
with a 3-core of size $q\le (3/2)\log{n}$, then 
by \cref{L_deg2s}
it has such a subgraph of size
$k=\beta\log{n}$ for some $\beta\in[\beta_*+\delta,1/\alpha]$. 
Select $\eta>0$ as in \cref{L_betaeps}. 
By \cref{L_eps*}, with high probability $\Gnp$
has no percolating $3$-core of size 
$q=\eps(2\alpha)^{-1}\log{n}$, 
for any $\eps\in[\eps_*+\eta,3\alpha]$.
On the other hand, by the choice of $\eta$,  
\cref{L_mueps,L_betaeps,L_mueps<mu}
imply that 
for any $\beta\in[\beta_*,1/\alpha]$, 
the expected number of 
irreducible percolating subgraphs of size $k=\beta\log{n}$
with a $3$-core of size $q\le(\eps_*+\eta)(2\alpha)^{-1}\log{n}$
is bounded by $n^{\mu+o(1)}$, where $\mu=\mu(\alpha,\beta)$. 
Hence the result follows by 
\cref{L_beta*}. 
 \end{proof}

\subsection{No percolating subgraphs 
with large cores}
\label{S_largeC}

To complete the proof of \cref{T_beta*}, 
we rule out 
the existence of large percolating $3$-cores.

\begin{prop}\label{P_3cores}
Fix $\alpha\in(0,1/3)$ and put $p=\p$. 
Then with high probability $\Gnp$ has no 
irreducible 
percolating $3$-cores $C$ 
of size $q=\beta\log{n}$, 
for any $\beta\in[3/2,9]$.
\end{prop}

Before proving the proposition we 
observe that it and \cref{P_smallC}
imply our main result. 

\begin{proof}[Proof of \cref{T_beta*}]
Let $\delta>0$ be given.  
By \cref{L_beta*}, 
we may assume that 
$\beta_*+\delta<3$.   
Hence, by \cref{L_k3k,L_core+deg2s}, 
if $\Gnp$ has a percolating subgraph 
that is larger 
than $(\beta_*+\delta)\log{n}$, then  
with high probability it 
has some irreducible percolating subgraph 
$G$ of size $k=\beta\log{n}$
with a core $C\subset G$ 
of size $q\le k$, for some 
$\beta\in(\beta_*+\delta,9]$.
By \cref{P_3cores}, with high probability 
$q\le (3/2)\log{n}$. However then, 
by \cref{P_smallC}, 
with high probability 
$\Gnp$ contains
no such subgraphs $G$. 
Therefore, with high probability,
all percolating subgraphs of $\Gnp$ are of size
$k\le(\beta_*+\delta)\log{n}$. 
On the other hand, as shown in \cite{AK18},  
$\Gnp$ has seed subgraphs of size 
larger than $(\beta_*-\delta)\log n$, completing the proof. 
 \end{proof}

Turning now to the proof of \cref{P_3cores}, we first 
observe that $\Gnp$ has no large percolating subgraphs  
with small cores and  
few vertices of degree 2.

\begin{lemma}\label{L_C3/2}
Fix $\alpha\in(0,1/3)$ and put $p=\p$. 
Fix some $i_*\ge1$. 
With high probability $\Gnp$ has no 
irreducible 
percolating subgraph $G$ 
of size $k\ge (3/2)\log{n}$
with a core $C\subset G$ 
of size $q\le (3/2)\log{n}$ and 
$i\le i_*$ vertices of degree $2$. 
\end{lemma}

This is essentially a straightforward 
consequence of \cref{L_Iellq}. 

\begin{proof}
By \cref{L_core+deg2s,L_extraE}, we may assume 
that if $\Gnp$ has an irreducible percolating subgraph $G$
of size $k=\beta\log{n}$ with a core of size $q\le(3/2)\log{n}$, 
then $G$ 
has excess $\ell\le 3$. 
By \cref{P_smallC,L_beta*,L_eps*}, 
we may further assume that 
$\beta\in[3/2,3]$ and  
$q=yk$, where $y\beta\in[0,3/2-\eps]$, 
for some $\eps>0$.
Without loss of generality, we assume that $\eps<1/(2e)$
and $\log(3/(2e)+\eps)<-1/2$
(which is possible, since $1+2\log(3/(2e))\approx-0.189<0$).
By \cref{L_Iellq} and since $\alpha<1/3$,  
for some constant $\vartheta\ge1$,
the expected number
of such subgraphs $G$ is bounded by 
\[
{n\choose k}p^{2k-3+\ell} I_q^\ell(k,i)
\le \vartheta k^{2\ell+i}p^{\ell-3} (knp^2\psi_\eps(q/k))^k
\ll n^{\nu}
\]
where
\[
\nu(\beta,\psi_\eps(y))=3/2+\beta\log(\beta/3)
+\beta\log\psi_\eps(y).
\]
Here, $\psi_\eps(y)$ is as defined in 
\cref{L_Iellq}, that is, 
\[
\psi_\eps(y)=
\max
\{
3/(2e)+\eps,(e/2)^{1-2y}y^2
\}.
\]
Recall that  
$\psi_\eps(y)=3/(2e)+\eps$
for $y\le y_*$ and $\psi_\eps(y)=(e/2)^{1-2y}y^2$
for $y>y_*$, where $y_*=y_*(\eps)$ is 
as defined by \eqref{E_y*}. Moreover,
$y_*\downarrow y_0$ as $\eps\downarrow0$, 
where $ y_0\approx0.819$. 

To complete the proof we show that, 
for some $\delta>0$,  
$\nu(\beta,\psi_\eps(y))<-\delta$  
for all relevant all $\beta,y$.
This follows by basic calculus, 
see \cref{A_C3/2}
below. 
 \end{proof}

Finally, we prove
\cref{P_3cores}.
The main idea is as follows: 
Suppose that $\Gnp$ has an irreducible  
percolating $3$-core $C$ of size 
$k=\beta\log{n}$, for some $\beta\in[3/2,9]$. 
By \cref{L_extraE}, we can assume that its excess is 
$\ell\le 27$. Hence, 
in the last step of a clique process for 
$C$, either 2 or 3 percolating subgraphs are merged
that have few vertices of degree 2
(by the observations following 
\cref{L_Iell} above). 
Therefore, by \cref{L_C3/2}, each of these subgraphs
is either smaller than $(3/2)\log{n}$,
or else has a $3$-core larger than 
$(3/2)\log{n}$. 
Hence, in proving \cref{P_3cores}, the key is consider  
$C$ as above of {\it minimal} size. 
By \cref{L_eps*}, there is some $\beta_1<3/2$ so that
with high probability $\Gnp$ has no percolating 
subgraphs
of size $\beta\log{n}$ with few vertices of degree 2, 
for $\beta\in[\beta_1,3/2]$. 
Hence such a graph $C$, if it exists, is the result of 
the (unlikely) event that 2 or 3 
percolating graphs, all of which are smaller than 
$\beta_1\log{n}$
and have few vertices of degree 2,   
are merged to form a 
percolating $3$-core larger than $(3/2)\log{n}$.
Informally, the existence of such a graph  
would require a ``macroscopic jump'' in 
the clique process. 

\begin{proof}[Proof of \cref{P_3cores}]
By \cref{L_eps*}, there is some $\beta_1<3/2$ so that
with high probability $\Gnp$ has no percolating 
subgraphs 
of size $\beta\log{n}$ with $i\le 15$ vertices of degree 2,
for $\beta\in[\beta_1,3/2]$.

Suppose that $\Gnp$ has an irreducible 3-core
$C$ of size $k=\beta\log{n}$, for some 
$\beta\in[3/2,9]$. 
By \cref{L_extraE}, we may assume that its
excess is $\ell\le 27$. 
Assume that $C$ is of minimal size amongst such 
subgraphs. 
Then by \cref{L_Emerge} there are two possibilities
for the last step of a clique process for $C$:
\begin{enumerate}[nolistsep,label=(\roman*)]
\item Three irreducible percolating 
subgraphs $G_j$, $j\in\{1,2,3\}$, 
are merged which form exactly one triangle 
$T=\{v_1,v_2,v_3\}$, such that for some 
$i_j\le 2$ and $k_j,\ell_j\ge0$ with $\sum k_j=k+3$
and $\sum\ell_j=\ell$,  
the 
$G_j$ contribute to $I^{\ell_j}(k_j,i_j)$.
If any $i_j>0$, the 
$i_j$  vertices of $G_j$ 
of degree 2 belong to $T$.

\item For some $m\le(\ell+3)/2\le15$, two 
percolating subgraphs $G_j$, $j\in\{1,2\}$,
are merged that share exactly $m$ vertices
$S=\{v_1,v_2,\ldots,v_m\}$, 
such that for some $i_j\le m$ and 
$k_j,\ell_j\ge0$  with $\sum k_j = k+m$ 
and $\sum\ell_j=\ell-(2m-3)$, 
the $G_j$ contribute to $I^{\ell_j}(k_j,i_j)$.  
If any $i_j>0$, the $i_j$
vertices of $G_j$ of degree 2 belong to $S$.  
\end{enumerate}

In either case, by the choice of $C$, all 
$G_j$ have a core smaller than $(3/2)\log{n}$.
Hence, by \cref{L_core+deg2s,L_extraE}, 
we may assume that 
each $\ell_j\le3$. Also, 
by \cref{L_C3/2}
and the choice of $\beta_1$,
we may further assume that all $G_j$ are 
smaller than 
$\beta_1\log{n}$.

{\bf Case (i)}.
Let $k,k_j,\ell_j$ be as in (i). 
Let $k_j-(j-1)=\eps_j k$, so that $\sum \eps_j=1$.
Without loss of generality 
we assume that $k_1\ge k_2\ge k_3$. 
Hence $\eps_1,\eps_2$ satisfy
$1/3\le\eps_1\le \beta_1/\beta<1$ and 
$(1-\eps_1)/2\le \eps_2\le\min\{\eps_1,1-\eps_1\}$. 
The number of $3$-cores $C$ as in (i) 
for these values $k,k_j,\ell_j$ is bounded by 
\[
{k\choose k_1,k_2-1,k_3-2}
{k_1\choose 2}{k_2-1\choose 1}2!^3
\prod_{j=1}^3
\sum_{i=0}^2{2\choose i}\frac{I^{\ell_j}(k_j,i)}{{k_j\choose i}}.
\] 
Applying \cref{L_Iell} and the inequality
$k!<ek(k/e)^k$ (and recalling $\ell_j\le3$), this is bounded by 
\[
{k\choose k-k_1}{k-k_1\choose k_3-2} 
\frac{k^3}{2}(8ek^7)^3
\left(\frac{2}{e^2}\right)^{k+3}
\prod_{j=1}^3 k_j^{2k_j}. 
\]
By the inequality ${n\choose k}<(ne/k)^k$,
and noting that 
\[
k_j^{2k_j}
\le
(ek)^{2(j-1)}
(k_j-(j-1))^{2(k_j-(j-1))},
\]
we see that the above expression is bounded by
$
(2e^{-2}\eta(\eps_1,\eps_2))^kk^{2k}n^{o(1)}
$, 
where
\begin{align*}
\eta(\eps_1,\eps_2)
&=
\left(\frac{e}{1-\eps_1}\right)^{1-\eps_1}
\left(\frac{(1-\eps_1)e}{\eps_3}\right)^{\eps_3}
\eps_1^{2\eps_1}
\eps_2^{2\eps_2}
\eps_3^{2\eps_3}\\
&=
\frac{e^{1-\eps_1+\eps_3}}{(1-\eps_1)^{\eps_2}}
\eps_1^{2\eps_1}\eps_2^{2\eps_2}\eps_3^{\eps_3}. 
\end{align*}
Since $\alpha<1/3$, 
the expected number of $3$-cores $C$ in $\Gnp$ of 
size $k=\beta\log{n}$ 
with $G_j$ of size $k_j$ as in (i)
is at most
\[
{n\choose k}p^{2k-3}\left(\frac{2}{e^2} 
\eta(\eps_1,\eps_2) k^2
\right)^k
n^{o(1)}
=p^{-3}\left(\frac{2}{e}\alpha\beta
\eta(\eps_1,\eps_2)\right)^kn^{o(1)}
\ll
n^{\nu}
\]
where
\[
\nu(\beta,\eps_1,\eps_2)
=\frac{3}{2}
+
\beta\log\left(\frac{2}{3e}\beta\eta(\eps_1,\eps_2)\right).
\]
Therefore, to show that with high probability
$\Gnp$ has no subgraphs $C$ as in (i) above, 
we need only show that, for some $\delta>0$, 
$\nu(\beta,\eps_1,\eps_2)<-\delta$
for all 
relevant $\beta,\eps_1,\eps_2$. 
This is proved in \cref{A_3cores}
by basic calculus. 

The next case is similar. We only sketch the details. 

{\bf Case (ii)}. 
Let $k,k_j,\ell_j,m$ be as in (ii). 
Let $k_1=\eps_1 k$ and $k_2-m=\eps_2 k$, 
so that $\sum \eps_j=1$.
Without loss of generality 
we assume that $k_1\ge k_2$. 
Hence  $\eps_1,\eps_2$ satisfy
$1/2\le\eps_1\le \beta_1/\beta<1$ and 
$\eps_2=1-\eps_1$. 
The number of 3-cores $C$ as in (ii) 
for these values $k,k_j,\ell_j,m$ is bounded by 
\[
{k\choose k_2-m}
{k_1\choose m}m!^2
\prod_{j=1}^2
\sum_{i=0}^m{m\choose i}\frac{I^{\ell_j}(k_j,i)}{{k_j\choose i}}.
\] 
Therefore, arguing as in Case (i), 
we find that 
the expected number of $3$-cores $C$ in $\Gnp$ of 
size $k=\beta\log{n}$ 
with $G_j$ of size $k_j$ as in (ii)
is $\ll n^\nu$, where
$\nu=\nu(\beta,\eps_1,1-\eps_1)$
is as in Case (i). 

The proof is complete.
 \end{proof}

\appendix

\section{Removing degree 2 vertices}\label{A_deg2s}

\begin{proof}[Proof of \cref{L_deg2s}]
The proof is by induction on the size of $G$.
The case $|V(G)|=3$, in which case $G$ is a triangle, is 
immediate.
Hence suppose that $G$, with $|V(G)|>3$, 
percolates and some $v\in V(G)$ is of degree 2, 
and assume that the statement of the lemma holds for all 
graphs $H$ with $|V(H)|<|V(G)|$. 

Let $(\CP_t)_{t=1}^\tau$ 
be a clique process for $G$. 
Let $e_1,e_2$ denote the edges incident to $v$ in $G$.
Let $t_v\le \tau$ be the first time in the clique process 
$(\CP_t)_{t=1}^\tau$
that a subgraph containing either $e_1$ or $e_2$ 
is merged with other (edge-disjoint, percolating) subgraphs. 
We claim that 
$\CP_{t_v}$ is obtained from $\CP_{t_v-1}$ 
by merging $e_1,e_2$ with
a subgraph in $\CP_{t_v-1}$.  
To see this, first note that if a graph $H$
percolates and $|V(H)|>2$ 
(i.e., $H$ is not simply an edge),
then all vertices in $H$ have degree at least 2.
Next, by the choice of $t_v$, observe that  
none of the graphs 
being merged contain both $e_1,e_2$. 
Therefore, since $v$ is of degree 2, 
if one the graphs contains 
an $e_i$,
it is necessarily equal to $e_i$. 
It follows that $v$ is contained in 
two of the graphs being merged,  
and hence that $\CP_{t_v}$ 
is the result of merging the edges $e_1,e_2$ with a subgraph 
in $\CP_{t_v-1}$, as claimed. 

To conclude, note that 
if $t_v=\tau$ then since $G$ percolates 
(and so $\CP_{t_v}=\{G\}$) 
we have that 
$\CP_{t_v-1}=\{e_1,e_2,G_ v\}$, 
and so $G_v$ percolates. 
Otherwise, if $t_v<\tau$, 
then $\CP_{\tau-1}$ 
consists of 2 or 3 subgraphs, 
one of which contains $e_1,e_2$. 
If $\CP_{\tau-1}=\{G_1,G_2\}$, 
where $e_1,e_2\in E(G_1)$, say, 
then 
$(G_1)_v$ 
percolates
by the inductive hypothesis. 
Since $G_1,G_2$ are edge-disjoint, 
we have that $v\notin V(G_2)$, 
as otherwise $G_2$
would be a percolating graph with an isolated vertex. 
Hence, by \cref{L_merge}(ii), 
$G_v=(G_1)_v\cup G_2$ percolates. 
Similarly, if $\CP_{\tau-1}=\{G_1,G_2,G_3\}$, 
where $e_1,e_2\in E(G_1)$, say, 
then by the inductive hypothesis 
and \cref{L_merge}(i), 
$G_v=(G_1)_v\cup G_2\cup G_3$ 
percolates. 
 \end{proof}

\section{Basic estimates}\label{A_Iell}

\begin{proof}[Proof of \cref{L_Iell}]
It is easily verified that the statement of the lemma
holds for $k\le 4$.  
For $k>4$, 
we claim moreover that for all $\ell\le3$ and relevant $i$, 
\begin{equation}\label{E_AIell}
I^\ell(k,i)\le A\zeta^k {k\choose i}k!k^{k+2\ell}
\end{equation}
where $\zeta=2/e$ and 
$A=6/(\zeta^55!5^5)$. 
Since 
$A<1$ and ${k\choose i}\le k^i$, 
the lemma follows. 

The constant $A$ 
is chosen as such 
to control the  case 
of 3-cores,  
$i=0$. 

The proof is by induction. 
By the choice of $A$, 
we note that 
\eqref{E_AIell}
holds for $k=5$. Indeed, 
$I(5,i)\le{5\choose i}{4\choose2}$ for all $i\in\{1,2,3\}$
and $I^\ell(5,i)=0$ otherwise. 
Assume that for some $k>5$, 
\eqref{E_AIell} holds for 
all $4<k'<k$, 
$\ell\le3$ and
relevant $i$. 

The case $i>0$, where $G$ has at least one vertex
of degree 2
follows easily, and 
explains the choice of $\zeta=2/e$.

{\bf Case 1} ($i>0$).  
Suppose that $G$ is a graph contributing to 
$I^\ell(k,i)$, where $i>0$ and $\ell\le3$. 
Let $v\in V(G)$ be the vertex of degree
2 in $G$ of minimal index. 
By considering which two of the $k-i$ vertices of $G$
of degree larger than 2
are neighbors of $v$, 
we find that $I^\ell(k,i)$ is bounded from above by 
\[
{k\choose i}{k-i\choose 2}
\sum_{j=0}^2
{2\choose j}
\frac{I^\ell(k-1,i-1+j)}{{k-1\choose i-1+j}}.
\]
In this sum, $j\in\{0,1,2\}$ is the number of neighbors of $v$ 
that are of degree 2 
in the subgraph of $G_v$ of $G$ induced by $V(G)-\{v\}$. 
Applying the inductive hypothesis,  we obtain 
\[
I^\ell(k,i)
\le 
A\zeta^k{k\choose i}k!k^{k+2\ell}
\cdot
\frac{2}{\zeta}\left(\frac{k-1}{k}\right)^{k}
\le A\zeta^k{k\choose i}k!k^{k+2\ell}, 
\]
as required. 

The remaining cases deal with $3$-cores $G$ of size $k$, 
where $i=0$. 
First, we establish the case $i=\ell=0$ of
edge-minimal $3$-cores.
The cases $i=0$ and $\ell\in\{1,2,3\}$
are proved by adapting this argument. 

{\bf Case 2} ($i=\ell=0$).  
Let $G$ be a graph contributing to $C(k)=I(k,0)$. 
Then, by \cref{L_Emerge}, 
in the last step of a clique process for $G$, 
three edge-minimal percolating subgraphs $G_j$,
$j\in\{1,2,3\}$, are merged which 
form exactly one triangle on some 
$T=\{v_1,v_2,v_3\}\subset V(G)$. Moreover, 
each $G_j$ has at most $2$ vertices of degree $2$, 
and if some $G_j$ has such a vertex $v$ then 
necessarily $v\in T$ (as else $G$ would have a vertex of degree 2). 
Also if $k_j=|V(G_j)|$, 
with $k_1\ge k_2\ge k_3$, 
then (i)
$\sum_{j=1}^3 k_j=k+3$, 
(ii) $k_1,k_2\ge4$ and  
(iii) $k_3=2$ or $k_3\ge 4$
(since if some $k_j=3$
or some $k_j=k_{j'}=2$, $j\neq j'$, 
then $G$ would have 
a vertex of degree 2). 

Since the inductive hypothesis only holds for graphs 
with more than 4 vertices, 
it is convenient to deal with the case $k_1=4$
separately: Note that the only 
irreducible percolating 3-cores of size $k$
with all $k_j\le4$ are of size $k\in\{7,9\}$. 
These graphs are the graph in \cref{F_3core}
and the graph obtained from this graph by 
replacing the bottom edge with a copy of $K_4$
minus an edge. It is easy 
to verify that \eqref{E_AIell}
holds if $k\in\{7,9\}$, 
and so in the arguments below 
we assume that $k_1>4$. 
Moreover, since the graph in \cref{F_3core} is the only
irreducible percolating 3-core on $k=7$ vertices, 
we further assume below that $k\ge8$.

We take three cases, with respect to 
whether (i) $k_2=4$, (ii) $k_2>4$ and $k_3\in\{2,4\}$, 
or (iii) $k_3>4$. 

{\bf Case 2(i)} ($i=\ell=0$ and $k_2=4$). 
Note that if $k_2=4$ then 
$k_3\in\{2,4\}$. 
The number of graphs $G$ as above with $k_3=2$
and $k_2=4$ is bounded from above by
\[
{k\choose k-3}{k-3\choose 2}{3\choose 1}2!^2
\sum_{j=0}^2
{2\choose j}
\frac{I(k-3,j)}{{k-3\choose j}}.
\]
Here the first binomial selects the vertices for the subgraph 
of size $k_1=k-3$, the next two binomials select the vertices
for the triangle $T$, and the rightmost factor bounds the number of 
possibilities for the subgraph of size $k_1=k-3$ 
(recalling that it can have at most 2 vertices of degree $2$, 
and if it contains any such vertex $v$, then $v\in T$). 
Applying the inductive hypothesis (recall that we may assume
that $k_1>4$), the above expression is bounded 
by 
\[
A\zeta^kk!k^k\cdot 
\frac{(k-3)^{k-1}}{k^k}\frac{4}{\zeta^3}
\le 
A\zeta^kk!k^k\cdot 
\frac{1}{k}\frac{4}{\zeta^3e^3}.
\]
Here, and throughout this proof, we use the fact that $(\frac{k-x}{k})^{k-y}\le e^{-x}$
provided that $2y\le x < k$ and $x>0$. To see this, note that 
$(\frac{k-x}{k})^{k-y}\to e^{-x}$ as $k\to\infty$, and 
\begin{align*}
\frac{\partial}{\partial k}\left(\frac{k-x}{k}\right)^{k-y}
&=\left(\frac{k-x}{k}\right)^{k-y}
\left(
\log\left(\frac{k-x}{k}\right)+\frac{x(k-y)}{k(k-x)}
\right)\\
&\ge\left(\frac{k-x}{k}\right)^{k-y}\frac{x(x-2y)}{2k(k-x)}\ge0,
\end{align*}
by the inequality $\log{u}\ge(u^2-1)/(2u)$ (which holds for $u\in(0,1]$). 

Similarly, the number of graphs $G$ as above such that  
$k_1=k_2=4$ is bounded by 
\[
{k\choose k-5,3,2}{k-5\choose 2}{3\choose 1}2!^3
\sum_{j=0}^2
{2\choose j}
\frac{I(k-5,j)}{{k-5\choose j}}. 
\]
By the inductive hypothesis, this is bounded by 
\[ 
A\zeta^kk!k^k\cdot 
\frac{(k-5)^{k-3}}{k^k}\frac{4}{\zeta^5}
\le
A\zeta^kk!k^k\cdot 
\frac{1}{k^{5/2}\sqrt{k-5}}\frac{4}{\zeta^5e^5}.
\]

Altogether, we find that the number of graphs $G$
contributing to $C(k)$ with $k_2=4$, divided by 
$A\zeta^kk!k^k$, is bounded by 
\begin{equation}\label{E_Gam1}
\gamma_1=
\frac{1}{8}\frac{4}{\zeta^3e^3}
+
\frac{1}{8^{5/2}\sqrt{3}}\frac{4}{\zeta^5e^5}
<0.07.
\end{equation}

{\bf Case 2(ii)} ($i=\ell=0$, $k_2>4$ and $k_3\in\{2,4\}$). 
Note that in this case we may further assume that $k\ge9$. 
For a given $k_1,k_2>4$,  the 
number of graphs $G$ as above
with $k_3=2$ (in which case $k_1+k_2=k+1$)
is bounded by
\[
{k\choose k_1,k_2-1}{k_1\choose 2}{k_2-1\choose 1}2!^2
\prod_{j=1}^2
\sum_{i=0}^2
{2\choose i}\frac{I(k_j,i)}{{k_j\choose i}}.
\]
Applying the inductive hypothesis, this is bounded 
by 
\[
A\zeta^kk!k^k\cdot 2\cdot4^2A\zeta\frac{k_1^{k_1+2}k_2^{k_2+2}}{k^k}.
\]
Since $k_2=k+1-k_1$, we have that 
\[
\frac{\partial}{\partial k_1}
k_1^{k_1+2}k_2^{k_2+2}
=-k_1^{k_1+1}k_2^{k_2+1}
(k_1k_2\log(k_2/k_1)-2(k_1-k_2)).
\]
By the bound $\log{x}\le x-1$, we see that 
\[
k_1k_2\log(k_2/k_1)-2(k_1-k_2)
\le -(k_2+2)(k_1-k_2)\le 0.
\]
Hence, setting $k_1$ to be the maximum relevant
value $k_1=k-4$ (when $k_2=5$), we find
\[
\frac{k_1^{k_1+2}k_2^{k_2+2}}{k^k}
\le 
\frac{5^7(k-4)^{k-2}}{k^k}\le \frac{1}{k^2}\frac{5^7}{e^4}
\]
for all relevant $k_1,k_2$. 
Therefore, summing over the at most $k/2$ possibilities
for $k_1,k_2$, we find that at most
\[
A\zeta^kk!k^k\cdot 
\frac{1}{k}\frac{A\zeta 4^2 5^7}{e^4}
\]
graphs $G$ with $k_3=2$ and $k_2>4$
contribute to $C(k)$. 

The case of $k_3=4$ is very similar. In this case, for a 
given $k_1,k_2>4$ such that $k_1+k_2=k-1$, 
the number of graphs $G$ as above is bounded by 
\[
{k\choose k_1,k_2-1,2}{k_1\choose 2}{k_2-1\choose 1}2!^3
\prod_{j=1}^2
\sum_{i=0}^2
{2\choose i}\frac{I(k_j,i)}{{k_j\choose i}},
\]
which, by the inductive hypothesis, is bounded 
by 
\[
A\zeta^kk!k^k\cdot 
2\cdot4^2\frac{A}{\zeta}\frac{k_1^{k_1+2}k_2^{k_2+2}}{k^k}.
\]
Arguing as in the previous case, we see that the above
expression is maximized when $k_2=5$ and $k_1=k-6$. 
Hence, summing over the at most $k/2$ possibilities
for $k_1,k_2$, there are at most
\[
A\zeta^kk!k^k\cdot
\frac{1}{(k-6)k^2}\frac{A 4^2 5^7}{\zeta e^6}
\]
graphs $G$ that contribute to $C(k)$ with $k_3=4$
and $k_2>4$. 

We conclude that the number of graphs $G$ that 
contribute to $C(k)$ with $k_2>4$
and $k_3\in\{2,4\}$, divided by 
$A\zeta^kk!k^k$, is bounded by 
\begin{equation}\label{E_Gam2}
\gamma_2=
\frac{1}{9}\frac{A\zeta 4^2 5^7}{e^4}
+
\frac{1}{3\cdot9^2}\frac{A 4^2 5^7}{\zeta e^6}
<0.15.
\end{equation}

{\bf Case 2(iii)} ($i=\ell=0$ and $k_3>4$). 
In this case we may further assume that $k\ge12$.
For a given $k_1,k_2,k_3>4$ such that $k_1+k_2+k_3=k+3$,  
the number of graphs $G$ as above
is bounded by
\[
{k\choose k_1,k_2-1,k_3-2}{k_1\choose 2}{k_2-1\choose 1}2!^3
\prod_{j=1}^3
\sum_{i=0}^2
{2\choose i}\frac{I(k_j,i)}{{k_j\choose i}}.
\]
By the inductive hypothesis, this is 
bounded by
\[
A\zeta^kk!k^k\cdot
2^24^3A^2\zeta^3
\frac{k_1^{k_1+2}k_2^{k_2+2}k_3^{k_3+2}}{k^k}.
\]
As in the previous cases considered, 
the above expression is maximized when 
$k_2=k_3=5$ and $k_1=k-7$. 
Hence, summing over the at most $k^2/12$ 
choices for the $k_j$, 
we find that at most 
\[
A\zeta^kk!k^k\cdot 
\frac{1}{((k-7)k)^{3/2}}\frac{A^2\zeta^3  4^3  5^{14}}{3e^7}
\]
graphs $G$ contribute to $C(k)$ with $k_3>4$. 
Hence, the number of such graphs, 
divided by $A\zeta^kk!k^k$, is bounded by 
\begin{equation}\label{E_Gam3}
\gamma_3=
\frac{1}{(5\cdot12)^{3/2}}\frac{A^2\zeta^3  4^3  5^{14}}{3e^7}
<0.01.
\end{equation}

Finally, combining 
\eqref{E_Gam1}, 
\eqref{E_Gam2}
and \eqref{E_Gam3},
we find that 
\begin{equation}\label{E_Gam123}
\frac{C(k)}{ A\zeta^kk!k^k}\le 
\gamma_1+\gamma_2+\gamma_3
< 0.23<1,
\end{equation}
completing the proof of Case 2. 

It remains to consider the cases $i=0$ 
and $\ell\in\{1,2,3\}$,
corresponding to $3$-cores $G$ with non-zero excess. 
In these cases, it is possible that 
only 2 subgraphs are merged
in the last step of a 
clique process for $G$.
We prove the cases $\ell=1,2,3$ separately, however
they all follow by adjusting the proof of Case 2.  

First, we note that  
if two graphs $G_1,G_2$ with at least 2 vertices in common
are merged to form 
an irreducible percolating 3-core $G$, then 
necessarily each $G_j$ contains more than 4
vertices. In particular, such a graph $G$ contains at least 8
vertices. 
This allow us to apply the inductive hypothesis
in these cases 
(recall that we claim that \eqref{E_AIell}
holds only for graphs with more than 4 vertices),
without taking additional sub-cases as in the proof 
of Case 2.

{\bf Case 3} ($i=0$ and $\ell=1$).
If $G$ contributes to $C^1(k)$, then 
by \cref{L_Emerge}, 
in the last step 
of a clique process for $G$, there are two cases
to consider: 
\begin{enumerate}[nolistsep,label={(\roman*)}]
\item Three percolating 
subgraphs $G_j$, $j\in\{1,2,3\}$, 
are merged which form exactly one triangle 
$T=\{v_1,v_2,v_3\}$, such that for some 
$i_j\le 2$ and $k_j,\ell_j\ge0$ with $\sum k_j=k+3$
and $\sum\ell_j=1$,  
we have that 
$G_j$ contributes to $I^{\ell_j}(k_j,i_j)$.
Moreover, if any $i_j>0$, the  $i_j$  vertices of $G_j$ 
of degree 2 belong to $T$.

\item Two percolating subgraphs $G_j$, $j\in\{1,2\}$,
are merged that share exactly two vertices
$S=\{v_1,v_2\}$, such that for some $i_j\le 2$ and 
$k_j$  with $\sum k_j = k+2$, we have that 
the $G_j$ contribute to $I(k_j,i_j)$. Moreover, 
if any $i_j>0$, the $i_j$
vertices of $G_j$ of degree 2 belong to $S$. 
\end{enumerate}

We claim that, by the arguments in Case~2 
leading to \eqref{E_Gam123}, 
the number of graphs $G$ satisfying (i),
divided by $A\zeta^kk!k^{k+2}$, 
is bounded by 
\begin{equation}\label{E_Gam123-ell1}
\gamma_1+2\gamma_2+3\gamma_3<0.40.
\end{equation}
To see this, note the only difference between (i) of the present case
and Case 2 above is that here one of the $G_j$ has
exactly 1 excess edge. 
Note that if one of the graphs $G_j$ 
has an excess
edge, then necessarily $k_j>4$.   
Recall that graphs 
$G$ that contribute to $C(k)$, as 
considered in Cases 2(i),(ii),(iii) above, 
have exactly $1,2,3$ subgraphs $G_j$ with $k_j>4$, 
respectively.
Moreover, recall that the number of such graphs $G$, 
divided by  $A\zeta^kk!k^k$, 
is bounded by $\gamma_1,\gamma_2,\gamma_3$, 
respectively, in these cases. 
Therefore, applying the inductive hypothesis, 
and noting that if $G_j$ 
has exactly $\ell_j=1$ excess edge then it 
contributes an extra factor of $k_j^2<k^2$, 
it follows that 
the number of graphs $G$ as in (i) of the present case, divided by 
$A\zeta^kk!k^{k+2}$, is bounded by $\sum_{j=1}^3j\gamma_j$,
as claimed. (By \eqref{E_Gam1}, \eqref{E_Gam2} and \eqref{E_Gam3}, 
this sum is bounded by $0.40$.)

On the other hand, arguing along the lines as in Case~2, 
the number of graphs $G$
satisfying (ii), for a given $k_1,k_2>4$ such that 
$k_1+k_2=k+2$, is bounded by 
\[
{k\choose k_1,k_2-2}{k_1\choose 2}2!^2
\prod_{j=1}^2
\sum_{i=0}^2{2\choose i}\frac{I(k_j,i)}{{k_j\choose i}}.
\]
By the inductive hypothesis, this is bounded by
\[
A\zeta^k k!k^k
\cdot 2\cdot 4^2A\zeta^2\frac{k_1^{k_1+2}k_2^{k_2+2}}{k^k}.
\] 
Arguing as in Case 2, we find that this expression is 
maximized when 
$k_2=5$ and $k_1=k-3$. Hence, 
summing over the at most $k/2$ choices for $k_1,k_2$, 
the number of graphs $G$ satisfying (ii),
divided by $A\zeta^kk!k^{k+2}$,
is at most
\begin{equation}\label{E_Gam4}
\gamma_4=\frac{1}{8^2}\frac{A\zeta^2 4^2 5^7}{e^3}
<0.04.
\end{equation}

Altogether, by \eqref{E_Gam123-ell1} and \eqref{E_Gam4}, 
we conclude that   
\begin{equation}\label{E_case3}
\frac{C^1(k)}{A\zeta^kk!k^{k+2}}
\le 
\gamma_1
+2\gamma_2
+3\gamma_3
+\gamma_4
<0.44<1,
\end{equation}
completing the proof of Case 3.

{\bf Case 4} ($i=0$ and $\ell=2$).
This case is nearly identical to Case 3. 
By \cref{L_Emerge}, in the last step of a 
clique process for a graph $G$
that contributes to $C^2(k)$, either (i) three graphs 
that form exactly one 
triangle are merged whose excesses sum to 2, 
or else (ii) two graphs that share exactly two 
vertices are merged whose excesses sum to 1. 
Hence, by the arguments
in Case 3 leading to \eqref{E_case3}, we find that 
\begin{equation}\label{E_case4}
\frac{C^2(k)}{A\zeta^k k!k^{k+4}}
\le 
\gamma_1
+3\gamma_2
+6\gamma_3
+2\gamma_4
< 0.66<1,
\end{equation}
as required.

{\bf Case 5} ($i=0$ and $\ell=3$).
Since $\ell=3$, it is now possible that 
in the last step of a clique process for a
graph $G$ contributing to $C^\ell(k)$, two 
graphs are merged that share three vertices. 
Apart from this difference, the argument is 
completely analogous to the previous cases. 

If $G$ contributes to $C^3(k)$, then by \cref{L_Emerge}, 
in the last step 
of a clique process for $G$, there are three cases
to consider: 
\begin{enumerate}[nolistsep,label={(\roman*)}]
\item Three percolating 
subgraphs $G_j$, $j\in\{1,2,3\}$, 
are merged which form exactly one triangle 
$T=\{v_1,v_2,v_3\}$, such that for some 
$i_j\le 2$ and $k_j,\ell_j\ge0$ with $\sum k_j=k+3$
and $\sum\ell_j=3$,  
we have that 
$G_j$ contributes to $I^{\ell_j}(k_j,i_j)$.
If any $i_j>0$, the 
corresponding $i_j$  vertices of $G_j$ 
of degree 2
belong to $T$.

\item Two percolating subgraphs $G_j$, 
$j\in\{1,2\}$,
are merged that share exactly two vertices
$S=\{v_1,v_2\}$, such that for some $i_j\le 2$ and 
$k_j,\ell_j\ge0$  with $\sum k_j = k+2$ and $\sum\ell_j=2$, 
we have that 
the $G_j$ contribute to $I^{\ell_j}(k_j,i_j)$.  
If any $i_j>0$, the $i_j$
vertices of $G_j$ of degree 2 belong to $S$. 

\item Two percolating subgraphs $G_j$, 
$j\in\{1,2\}$,
are merged that share exactly three vertices
$R=\{v_1,v_2,v_3\}$, such that for some $i_j\le 3$ 
and $k_j$ 
with $\sum k_j = k+3$, 
we have that 
the $G_j$ contribute to $I(k_j,i_j)$.  
If any $i_j>0$, the $i_j$
vertices of $G_j$ of degree 2 belong to $R$.
\end{enumerate}

As in Case 4, we find by the arguments
in Case 3 leading to \eqref{E_case3} that the number of 
graphs $G$ satisfying (i) or (ii),
divided by $A\zeta^k k!k^{k+6}$, is bounded by 
\begin{equation}\label{E_case5i}
\gamma_1
+4\gamma_2
+10\gamma_3
+3\gamma_4
<0.89. 
\end{equation}

By the arugments in Case 3 leading to 
\eqref{E_Gam4}, 
the number of graphs $G$
satisfying (iii), for a given $k_1,k_2>4$ such that 
$k_1+k_2=k+3$, is bounded by
\[
{k\choose k_1,k_2-3}{k_1\choose 3}3!^2
\prod_{j=1}^2
\sum_{i=0}^3
{3\choose i}\frac{I(k_j,i)}{{k_j\choose i}}.
\]
By the inductive hypothesis, this is bounded by
\[
A\zeta^k k!k^k\cdot  3! 8^2
A\zeta^3\frac{k_1^{k_1+3}k_2^{k_2+3}}{k^k}.
\] 
This expression is 
maximized when 
$k_2=5$ and $k_1=k-2$. Hence, 
summing over the at most $k/2$ choices for $k_1,k_2$, 
the number of graphs $G$ satisfying (iii),
divided by $A\zeta^kk!k^{k+6}$,
is at most
\begin{equation}\label{E_Gam5}
\gamma_5
=\frac{1}{8^4}\frac{A\zeta^3 3! 5^8 8^2}{2e^2}
<0.08.
\end{equation}

Therefore, by \eqref{E_case5i} and \eqref{E_Gam5},
we have that 
\[
\frac{C^3(k)}{A\zeta^k k!k^{k+6}}
\le 
\gamma_1
+4\gamma_2
+10\gamma_3
+3\gamma_4
+\gamma_5
< 0.97<1,
\]
completing the proof of Case 5.

This last case completes the induction.  
We conclude that \eqref{E_AIell}
holds for all $k>4$, $\ell\le 3$ and relevant $i$,
and the lemma follows. 
 \end{proof}

\section{Sharper estimates}\label{A_Iellq}

\begin{proof}[Proof of \cref{L_Iellq}]
Let $\eps>0$ be given. 
We may assume that $\eps<1/(2e)$, as otherwise
the statement of lemma follows by \cref{L_Iell}.
We claim that, for some $\vartheta(\eps)\ge1$ (to be determined below), 
and for all $k\ge2$, $\ell\le 3$ and relevant $q,i$, we have that 
\begin{equation}\label{E_psi}
I_q^\ell(k,i)\le \vartheta {k\choose i}\psi_\eps(q/k)^k k!k^{k+2\ell}.
\end{equation}

{\bf Case 1} ($i=k-q$). 
We first observe that \cref{L_Iell}
implies the case $i=k-q$. 
Indeed, if $q=k$, in which case $i=0$, then \eqref{E_psi} 
follows immediately by \cref{L_Iell}, 
noting that $I_k^\ell(k,0)=C^\ell(k)$ and $\psi(1)=2/e$.
On the other hand, if $i=k-q>0$ then 
\[
I^\ell_q(k,k-q) 
= {k\choose k-q}{q\choose 2}^{k-q}C^\ell(q),
\]
since all $k-q$ vertices of degree 2 in a graph 
that contributes to $I^\ell_q(k,k-q)$ are neighbors of 2 vertices
in its core. 
We claim that the right hand side is bounded by 
\[
{k\choose k-q}
(e/2)^{k-2q}(q/k)^{2k}
k!k^{k+2\ell}.
\]
Since $(e/2)^{k-2q}(q/k)^{2k}\le \psi(q/k)^k$,  
\eqref{E_psi} follows. 
To see this, note that by \cref{L_Iell},  
we have that 
\[
\frac{{q\choose 2}^{k-q}C^\ell(q)}{(e/2)^{k-2q}(q/k)^{2k}k!k^{k+2\ell}}
\le
\left(\frac{q}{k}\right)^{2\ell}
\frac{q!}{(q/e)^q}\frac{(k/e)^k}{k!}
\le 
\frac{q!}{(q/e)^q}\frac{(k/e)^k}{k!}.
\]
By the inequalities
$1\le i!/(\sqrt{2\pi i}(i/e)^i)\le e^{1/(12i)}$, it is easy to verify 
that the right hand side above is bounded by 1, for all relevant $q\le k$. 
Hence \eqref{E_psi} holds also in the case $i=k-q>0$. 

{\bf Case 2} ($i<k-q$). 
Fix some $k_\eps\ge 1/(1-y_*)^2$ 
(where $y_*$ is as in \eqref{E_y*}) such that, for all $k\ge k_\eps$
and relevant $q$, we have that 
\[
1+\frac{2}{k-1}\left(\frac{k-2}{k-1}\right)^k
\frac{\psi_\eps(q/(k-2)^{k-2}}{\psi_\eps(q/(k-1))^{k-1}}
=1+O(1/k)\le 1+\delta,
\]
where 
\[
\delta=\min
\left\{
1-\frac{3/(2e)}{3/(2e)+\eps},1-\frac{3(1-y_*)}{y_*^2}
\right\}.
\]
Note that, since $3(1-y)/y^2<1$ for all
$y>(\sqrt{21}-3)/2\approx 0.791$, 
and recalling (see \eqref{E_y*}) 
that $y_*> y_0\approx0.819$, it follows that  $\delta>0$. 

Select $\vartheta\ge1$ 
so that \eqref{E_psi} holds for all $k\le k_\eps$
and relevant $q,\ell,i$. 
By Case 1 and since $\vartheta\ge1$, we have that 
\eqref{E_psi} holds for all $k,q$ in the case that $i=k-q$. 
We establish the remaining cases $i<k-q$ by induction. 
Assume that for some $k>k_\eps$, 
\eqref{E_psi} holds for all $k'<k$ 
and relevant $q,\ell,i$. 

In any graph $G$ contributing to $I_q^\ell(k,i)$, 
where $i<k-q$, 
there is some vertex
of degree 2 with at least one of its two neighbors 
not in the core of $G$. 
There are two cases to consider: either
\begin{enumerate}[nolistsep,label=(\roman*)]
\item there is a vertex $v$
of degree 2 such that at least one of its two neighbors 
is of degree 2 in $G_v$ (obtained from $G$ by deleting $v$), 
or else,  
\item there is no such vertex $v$, 
however there are vertices $u\neq w$
of degree 2 in $G$ with a common neighbor that is not in the
core $C$ of $G$. 
\end{enumerate}
Note that, in case (i), removing $v$ 
results in a graph with $j\in\{i,i+1\}$ vertices of degree 2. 
On the other hand, 
in case (ii), removing $u$ and $w$ results in a graph with $j\in\{i-2,i-1,i\}$
vertices of degree 2. 
Hence, for $i<k-q$, 
we find that 
$I^\ell_q(k,i)/{k\choose i}$ is bounded by 
\begin{multline*}
\frac{I^\ell_q(k-1,i+1)}{{k-1\choose i+1}}{k-i-q\choose2}
+
\frac{I^\ell_q(k-1,i)}{{k-1\choose i}}(k-i-q)(k-i)\\
+
(k-i-q)(k-i)^2
\sum_{j=0}^2
\frac{I^\ell_q(k-2,i-2+j)}{{k-2\choose i-2+j}}.
\end{multline*} 
Applying the inductive hypothesis, it follows 
(after simple, but somewhat tedious simplifications) 
that 
\[
\frac{I_q^\ell(k,i)}{\vartheta {k\choose i}\psi_\eps(q/k)^k k!k^{k+2\ell}}
\le
\Psi_\eps(q,k)
\left[
1+\frac{2}{k-1}\left(\frac{k-2}{k-1}\right)^k
\frac{\psi_\eps(q/(k-2)^{k-2}}{\psi_\eps(q/(k-1))^{k-1}}
\right]
\]
where
\[
\Psi_\eps(q,k)=
\frac{3}{2}\frac{k-q}{k}\left(\frac{k-1}{k}\right)^k
\frac{\psi_\eps(q/(k-1))^{k-1}}{\psi_\eps(q/k)^k}. 
\]
By the choice of $k_\eps$, and since $k\ge k_\eps$,
we have that 
\begin{equation}\label{E_psi-ind}
\frac{I_q^\ell(k,i)}{\vartheta{k\choose i} \psi_\eps(q/k)^k k!k^{k+2\ell}}
\le
\Psi_\eps(q,k)
(1+\delta).
\end{equation}

Next, we show that $\Psi_\eps(q,k)<1-\delta$,
completing the induction.
To this end, we take cases with respect to whether 
(i) $q/(k-1)\le y_*$, (ii) $y_*\le q/k$,  or 
(iii) $q/k<y_*<q/(k-1)$.  

{\bf Case 2(i)} ($q/(k-1)\le y_*$).
In this case $\psi_\eps(q/m)=3/(2e)+\eps$, 
for each $m\in\{k-1,k\}$. 
It follows, by the choice of $\delta$, that 
\[
\Psi_\eps(q,k)
\le 
\left(\frac{k-1}{k}\right)^k\frac{3/2}{3/(2e)+\eps}
\le
\frac{3/(2e)}{3/(2e)+\eps}
<1-\delta, 
\]
as required. 

{\bf Case 2(ii)} ($y_*\le q/k$).
In this case, 
we have that $\psi(q/m)^m=(e/2)^{m-2q}(q/m)^{2m}$,
for each $m\in\{k-1,k\}$.
Hence
\[
\Psi_\eps(q,k)
=
\frac{3}{e}\left(\frac{k}{k-1}\right)^{k-1}\frac{(k-q)(k-1)}{q^2}
\le \frac{3(1-y)}{y^2},  
\]
where $y=q/k$.
Since the right hand side is decreasing in $y$, we find, by the choice
of $\delta$, that 
\[
\Psi_\eps(q,k)\le \frac{3(1-y_*)}{y_*^2}
<1-\delta.
\]

{\bf Case 2(iii)} ($q/k<y_*<q/(k-1)$).   
In this case, $\psi_\eps(q/k)=3/(2e)+\eps$ and 
\[
\psi_\eps(q/(k-1))^{k-1}=(e/2)^{k-1-2q}(q/(k-1))^{2(k-1)}.
\]
Hence
\[
\Psi_\eps(q,k)
=
\frac{3}{e}\left(\frac{k}{k-1}\right)^{k-1}\frac{(k-q)(k-1)}{q^2}
\frac{(e/2)^{k-2q}(q/k)^{2k}}{(3/(2e)+\eps)^k}. 
\]
As in the previous case, we consider the quantity $y=q/k$. 
The above expression is bounded by 
\[
\frac{3(1-y)}{y^2}\left(
\frac{(e/2)^{1-2y}y^2}{3/(2e)+\eps}
\right)^k.
\]
We claim that this expression is increasing in $y\le y_*$.
By \eqref{E_y*} and the choice of $\delta$, 
it follows that 
\[
\Psi_\eps(q,k)\le \frac{3(1-y_*)}{y_*^2}
<1-\delta,
\]
as required.
To establish the claim, 
simply note that 
\begin{align*}
\frac{\partial}{\partial y}
\frac{1-y}{y^2}
((2/e)^y y)^{2k}
&=
\frac{1}{y^3}((2/e)^y y)^{2k}
\left(
2(1-y)(1+y\log(2/e))k+y-2
\right)\\
&>
\frac{2}{y^3}((2/e)^y y)^{2k}
(  (1-y)^2k-1  )\ge0
\end{align*}
for all $y\le y_*$, since $k\ge k_\eps\ge 1/(1-y_*)^2$.

Altogether, we conclude that 
$\Psi_\eps(q,k)\le 1-\delta$, for all relevant $q$. 
By \eqref{E_psi-ind}, it follows that 
\[
\frac{I_q^\ell(k,i)}{\vartheta {k\choose i}\psi_\eps(q/k)^k k!k^{k+2\ell}}
\le
1-\delta^2
<
1
\]
completing the induction. We conclude that \eqref{E_psi}
holds for $k\ge2$, $\ell\le 3$ and relevant $q,i$. 
Since ${k\choose i}\le k^i$, the lemma follows. 
 \end{proof}

\section{Details in 
the proof of \cref{L_C3/2}}\label{A_C3/2}

In this section, to complete the proof of 
\cref{L_C3/2}, we verify that, for some $\delta>0$, 
we have that 
$\nu(\beta,\psi_\eps(y))<-\delta$  
for all relevant $\beta,y$. 
Note that $\nu$ is 
convex in $\beta$. 
Therefore it suffices to consider the extreme points 
$\beta=3/2$ and $\beta=\min\{3,3/(2y)\}$
in the range $y\in[0,1-2\eps/3]$. 

Since $\psi_\eps(1)=2/e$, we have that  
$\nu(3/2,\psi_\eps(1))=0$.
Hence, for some $\delta_1>0$, 
we have that 
$\nu(3/2,\psi_\eps(y))<-\delta_1$ for all
$y\in[0,1-2\eps/3]$. 
Next, 
for $\beta=\min\{3,3/(2y)\}$, we treat
the cases (i) $y\in[0,1/2]$ and $\beta=3$ and 
(ii) $y\in[1/2,1-2\eps/3]$ and $\beta=3/(2y)$ separately. 
If $y\le1/2$, then $\psi_\eps(y)=3/(2e)+\eps$,
in which case, by the choice of $\eps$, 
\[
\nu(3,\psi_\eps(y))
=\frac{3}{2}(1+2\log(3/(2e)+\eps))
<0. 
\]
On the other hand, for $y\ge1/2$, 
we need to show that  
\[
\nu(3/(2y),\psi_\eps(y))=\frac{3}{2}\left(
1+\frac{1}{y}\log\left(\frac{\psi_\eps(y)}{2y}\right)
\right)<0.
\]
To this end, we first note that differentiating 
$\nu(3/(2y),3/(2e)+\eps)$
twice with respect to $y$, we obtain 
\[
\frac{3}{2y^3}\left(3
+2\log\left(\frac{3/(2e)+\eps}{2y}
\right)
\right)
\ge
\frac{3}{2}\left(3
+2\log\left(\frac{3}{4e}
\right)
\right)
\approx 0.637
>0.
\]
Therefore it suffices to consider the extreme points 
$y=1/2$ and $y=1$. Noting that,
by the choice of $\eps$, 
we have that 
\[
\nu(3,3/(2e)+\eps)
=\frac{3}{2}(1+2\log(3/(2e)+\eps))
<0
\]
and
\begin{align*}
\nu(3/2,3/(2e)+\eps)
&=\frac{3}{2}\left(1+\log\left(\frac{3/(2e)+\eps}{2}\right)\right)\\
&<\frac{3}{2}(1+2\log(3/(2e)+\eps))
<0, 
\end{align*}
it follows that $\nu(3/(2y),3/(2e)+\eps)<0$ for all $y\in[1/2,1]$. 
Next, we observe that differentiating 
$\nu(3/(2y),(e/2)^{1-2y}y^2)$ with respect to $y$, 
we obtain 
\[
\frac{3}{2y^2}
\left(
1-\log(ey/4)
\right)
\ge
3\log{2}
>0.
\]
Therefore, since 
$\nu(3/(2y),(e/2)^{1-2y}y^2)\to\nu(3/2,\psi_\eps(1))=0$
as $y\uparrow1$, it follows that 
$\nu(3/(2y),(e/2)^{1-2y}y^2)<0$ for all $y\in[1/2,1-2\eps/3]$. 
Altogether, there is some 
$\delta_2>0$ so that 
$\nu(\min\{3,3/(2y)\},\psi_\eps(y))<-\delta_2$
for all $y\in[0,1-2\eps/3]$. 

Taking $\delta=\min\{\delta_1,\delta_2\}$, 
it follows that 
$\nu(\beta,\psi_\eps(y))<-\delta$, 
for all relevant $\beta,y$, as required. 

\section{Details in the proof of 
\cref{P_3cores}}\label{A_3cores}

We finish the proof of \cref{P_3cores}
by showing that, for some $\delta>0$, we have 
$\nu(\beta,\eps_1,\eps_2)<-\delta$, 
for all relevant $\beta,\eps_1,\eps_2$. 
Since $\nu$ is
convex in $\beta$, we can restrict to 
the extreme points $\beta=3/2$
and $\beta=3/(2\eps_1)>\beta_1/\eps_1$.
To this end, observe that when $\beta=3/2$, 
we have that $\nu<0$ if and only if
$\eta<1$. Similarly, when $\beta=3/(2\eps_1)$, 
$\nu<0$ if and only if
$\eta<\eps_1e^{1-\eps_1}$. Since $\eps_1e^{1-\eps_1}\le1$
for all relevant $\eps_1$, it suffices to establish the latter claim.
To this end, we observe that  
\begin{align*}
\frac{\partial}{\partial \eps_2} \eta(\eps_1,\eps_2)
&=\eta(\eps_1,\eps_2)
\log\left(
\frac{e\eps_2^2}{(1-\eps_1)(1-\eps_1-\eps_2)}
\right)\\
&\ge
\eta(\eps_1,\eps_2)\log(e/2)
>0
\end{align*}
for all relevant $\eps_2\ge(1-\eps_1)/2$. 
Therefore, we need only show that 
\[
\zeta(\eps_1)
=
\frac{\eta(\eps_1,\min\{\eps_1,1-\eps_1\})}{\eps_1e^{1-\eps_1}}
<1-\delta
\]
for some $\delta>0$ and all relevant $\eps_1$.
We treat the cases $\eps_1\in[1/3,1/2]$ and 
$\eps_1\in[1/2,1)$ separately. 

For $\eps_1\in[1/3,1/2]$, we have
\[
\zeta(\eps_1)
=
\frac{\eta(\eps_1,\eps_1)}{\eps_1e^{1-\eps_1}}
= 
\frac{(e(1-2\eps_1))^{1-2\eps_1}\eps_1^{4\eps_1-1}}{(1-\eps_1)^{\eps_1}}.
\]
Hence 
\[
\frac{\partial}{\partial\eps_1}\zeta(\eps_1)
=\zeta(\eps_1)
\left(\log\left(\frac{\eps_1^4}{(1-\eps_1)(1-2\eps_1)^2}\right)+
\frac{\eps_1^2+\eps_1-1}{\eps_1(1-\eps_1)}\right).
\]
The terms $\eps_1^4/((1-\eps_1)(1-2\eps_1)^2)$
and $(\eps_1^2+\eps_1-1)/(\eps_1(1-\eps_1))$
are increasing for $\eps_1\in[1/3,1/2]$, as is easily verified. 
Hence $\zeta(\eps_1)$ is decreasing in $\eps_1$
for $1/3\le\eps_1\le x_1\approx 0.439$ and 
increasing for $x_1\le \eps_1\le1/2$. 
Therefore, since  $\zeta(1/3)=(e/6)^{1/3}<1$
and $\zeta(1/2)=1/\sqrt{2}<1$, 
we have that, for some $\delta_1>0$, 
$\zeta(\eps_1)<1-\delta_1$
for all $\eps_1\in[1/3,1/2]$.

Similarly, for $\eps\in[1/2,1)$, we have
\[
\zeta(\eps_1)
=
\frac{\eta(\eps_1,1-\eps_1)}{\eps_1e^{1-\eps_1}}
=
(1-\eps_1)^{1-\eps_1}\eps_1^{2\eps_1-1}. 
\]
Hence 
\[
\frac{\partial}{\partial\eps_1}\zeta(\eps_1)
=\zeta(\eps_1)
\left(\log\left(\frac{\eps_1^2}{1-\eps_1}\right)
+\frac{\eps_1-1}{\eps_1}\right). 
\]
Since $\eps_1^2/(1-\eps_1)$ and $(\eps_1-1)/\eps_1$
are increasing in $\eps_1\in[1/2,1)$, 
we find that $\zeta(\eps_1)$ is decreasing in $\eps_1$
for $1/2\le \eps_1\le x_2\approx 0.692$ and 
increasing for $x_2\le \eps_1<1$. 
Note that 
$\zeta(1/2)=1/\sqrt{2}<1$
and 
$\zeta(1)=1$. 
Hence, for some $\delta_2>0$, 
$\zeta(\eps_1)<1-\delta_2$ for all
$\eps_1\in[1/2,\beta_1/\beta]\subset[1/2,1)$. 

Setting $\delta'=\min\{\delta_1,\delta_2\}$, 
we find that 
$\zeta(\eps_1)<1-\delta'$
for all relevant $\eps_1$. 
It follows that, for some $\delta>0$, 
we have that 
$\nu(\beta,\eps_1,\eps_2)<-\delta$, 
for all relevant $\beta,\eps_1,\eps_2$.

\providecommand{\bysame}{\leavevmode\hbox to3em{\hrulefill}\thinspace}
\providecommand{\MR}{\relax\ifhmode\unskip\space\fi MR }
\providecommand{\MRhref}[2]{%
  \href{http://www.ams.org/mathscinet-getitem?mr=#1}{#2}
}
\providecommand{\href}[2]{#2}

\end{document}